\documentclass[journal]{IEEEtran}
\usepackage{cite}
\usepackage{url}
\usepackage{amsmath,amssymb,amsfonts}
\usepackage{algorithmic}
\usepackage{graphicx}
\usepackage{textcomp}
\usepackage{mathtools}
\usepackage{subfigure}
\usepackage[misc]{ifsym}
\usepackage{color}
\usepackage[numbers,sort&compress]{natbib}

\def\BibTeX{{\rm B\kern-.05em{\sc i\kern-.025em b}\kern-.08em
    T\kern-.1667em\lower.7ex\hbox{E}\kern-.125emX}}
\begin{document}
\title{Market Power in Convex Hull Pricing}
\author{Jian Sun \emph{and} Chenye Wu\vspace{-0.8cm}
\thanks{The authors are with the Institute of Interdisciplinary Information Sciences (IIIS), Tsinghua University, Beijing, China, 100084. C. Wu is the correspondence author. Email:chenyewu@tsinghua.edu.cn}}

\maketitle

\begin{abstract}
The start up costs in many kinds of generators lead to complex cost structures, which in turn yield severe market loopholes in the locational marginal price (LMP) scheme. Convex hull pricing (a.k.a. extended LMP) is proposed to improve the market efficiency by providing the minimal uplift payment to the generators. In this letter, we consider a stylized model where all generators share the same generation capacity. We analyze the generators' possible strategic behaviors in such a setting, and then propose an index for market power quantification in the convex hull pricing schemes.
\end{abstract}

\section{Introduction}

In general, the electricity market is organized by a sequence of processes: unit commitment, day-ahead market, real time market, balancing market, \emph{etc}. Since the focuses of different processes are diverse, the incentives they provide do not necessarily align with each other. To relieve the tension between these processes, convex hull pricing (CHP) is one promising solution: it provides the minimal uplift payment to the generators for incentive alignment \cite{o2005efficient,hogan2003minimum,Schiro}.

However, although remarkable in improving the system efficiency by reallocating the market surplus between supply and demand \cite{Thompson}, we find that the incentive issues in the CHP scheme have not been fully solved due to the non-convex cost structures. To the best of our knowledge, we are the first to identify market power in the CHP schemes.

In this letter, we first revisit the CHP scheme in Section \ref{sec:2}. And then we seek to understand each individual firm's possible strategic behaviors in Section \ref{sec:3}. Based on the identified behaviors, we propose the market power index and highlight the existence of market power via numerical studies in Section \ref{sec:4}. Finally, concluding remarks are given in Section \ref{sec:5}. We provide all the necessary proofs in the Appendix.

\section{Convex Hull Pricing: A Brief Revisit}
\label{sec:2}
We consider the CHP scheme in the electricity pool model with $n$ generators in the system. To highlight the existence of strategic behaviors, we assume all the generators have the same capacity $G$. In such a model, the system operator seeks to minimize the total generation costs at each time slot:
\begin{equation}
    \begin{aligned}
    c(y)\coloneqq \min_\textbf{g} &\sum\nolimits_{i=1}^n f_i(g_i)\\
    s.t. & \sum\nolimits_{i=1}^n g_i=y,\\
    &0\le g_i\le G,
\end{aligned}
\end{equation}
where $g_i$ denotes the output of generator $i$, $f_i(g_i)$ denotes the generation cost for generator $i$, vector $\textbf{g}$ is $[g_1,...,g_n]$, and $y$ is the total load in the system. The optimal solution $\textbf{g}^*(y)$ and the optimal objective value $c(y)$ are both functions of $y$.

If $f_i(g_i)$ is linear, then the operator could conduct the economic dispatch according to the merit order of marginal generation cost and set the price as the marginal cost.

However, in practice, there are start up costs associated with generation. Hence, for generator $i$, $f_i(g_i)$ is often of the following form:
\begin{equation}
    f_i(g_i)=s_i\cdot I(g_i>0)+v_i\cdot g_i,
    \label{eq:cost structure}
\end{equation}
where $s_i$ is the sum of fixed cost and start up cost, and $v_i$ is the variable cost. To highlight the non-convexity in $f_i(g_i)$, we employ the indicator function $I(\cdot)$. Such a cost structure challenges the conventional scheme in terms of dispatch profile, price design, and incentive analysis.

The solution proposed by CHP is to use the uplift payment to incentivize the generators to follow the system operator's dispatch profile.

For each generator $i$, given price $p$, its desired generation level $g_i^d(p)$ is associated with its maximal profits, and can be different from the dispatch profile $g_i^*$:
\begin{equation}
    g_i^d(p)=\sup \left\{\arg \max_{z\in [0,G]}\{p\cdot z-f_i(z)\}\right\}.
\end{equation}
And this generation level leads to the maximal profits given $p$, denoted by $\pi_i^*(p)$, i.e.,
 \begin{equation}
     \pi_i^*(p)=p\cdot g_i^d(p)-f_i(g_i^d(p)).
 \end{equation}
 The difference in the profits between generating $g_i^d(p)$ and $g_i^*$ needs to be compensated by uplift payment:
\begin{equation}
    \text{Uplift}_i(p,y)=\pi_i^*(p)-[p\cdot g_i^*-f_i(g_i^*)].
    \label{def:uplift}
\end{equation}
And the CHP scheme is proved to guarantee that it leads to the optimal price $p^*$ in terms of minimizing the total uplift payment \cite{Gribik}. In this letter, we adopt an equivalent form to characterize $p^*$:

\vspace{0.1cm}
\noindent \textbf{Lemma 1}: The optimal price $p^*$ in CHP can be expressed as follows:
\begin{equation}
    p^*=\inf\left\{p\big |\sum\nolimits_{i=1}^n g_i^d(p)\ge y\right\}.
\end{equation}
Given the structure of $f_i(g_i)$, we claim for demand of $y$, there exists an optimal dispatch profile $\textbf{g}^*$, which is composed of the following dispatch profiles:

1) $m-1$ elements are $G$, and $n-m$ elements are 0.

2) one generator is dispatched at $x$ units. 

Here, $m=\lceil \frac{y}{G}\rceil$, and $x=y-(m-1)\cdot G$. This allows the system operator to dispatch the generators roughly in the order of their average generation costs\footnote{This order plays the most important role in the optimal dispatch. There is minor issue in dispatching the generator which is scheduled to generate $x$ units. We omit the detailed discussion due to the page limits.}. In our case, since all the generators share the same capacity, the system operator can directly sort them with respect to $f_i(G)$. For notational simplicity, we assume the subscript $i$ for generator $i$ also denotes its ascending order in the average generation costs. Hence, we can make the following observations.

\vspace{0.1cm}
\noindent \textbf{Fact 1}: We can establish the mapping between generator's order in average generation cost and its dispatch profile.

1) If $i<m$, then $g_i^*\in\{x,G\}$; if $i>m$, then $g_i^*\in\{0,x\}$.

2) If $g_i^*=0$, then $i\ge m$; if $g_i^*=G$, then $i\le m$.
\vspace{0.1cm}

However, the CHP scheme only guarantees the minimal total uplift payment, it does not fully mitigate generators' opportunities in obtaining more profits by strategic bidding.

\section{Strategic Behavior Analysis}
\label{sec:3}
 We seek to understand generator's strategic behavior via analyzing its profits. For each generator, its total profits are gained from the uplift payment as well as selling electricity according to the dispatch profile from the system operator. We want to emphasize that the second component is not necessarily positive, which further implies the importance of uplift payment. 
 
 In practice, the generators are allowed to bid their cost functions as well as their available capacities to the operator. In this letter, we assume the generator is not allowed to withhold its capacity. Hence, it can only strategically bid its cost function. Since the fixed cost and start up cost are rather stable, the only remaining opportunity for manipulation is the variable cost $v_i$, for each generator $i$.
 
 Mathematically, if generator $i$ truthfully bids its $v_i$, given a demand of $y$, we denote its profits $P_i$ as benchmark:
 \begin{equation}
     P_i = \pi_i(p^*(v_i)) = \{f_m(G)-f_i(G)\}^+.
     \label{eq:truthful profit}
 \end{equation}

 Recall $m=\lceil \frac{y}{G}\rceil$, and $f_m(\cdot)$ denotes the cost function for generator $m$, which ranks the $m^{th}$ in terms of the average generation cost. It is crucial in determining the optimal price $p^*$ in CHP. We provide the detailed analysis in the Appendix.
 
 On the other hand, if the generator strategically reports its generation cost as $\Tilde{f}_i(\cdot)$ (more precisely, a different variable cost $\tilde{v}_i$), which may lead to a potentially different $\Tilde{p}^*$ given by CHP, and a potentially different dispatch profile $\Tilde{\textbf{g}}^*$. This strategic bidding will also reshuffle the order of generators in terms of average \emph{bid} generation cost. We denote the reshuffled bid generation cost by $\tilde{f}^{(1)},\cdots,\tilde{f}^{(n)}$. Hence, generator $i$'s total profits $\Tilde{P}_i(\tilde{v}_i)$ via strategic bidding by can be straightforwardly characterized as follows:
\begin{equation}
    \Tilde{P}_i(\tilde{v}_i)=\underbrace{\left\{\tilde{f}^{(m)}(G)-\tilde{f}_i(G)\right\}^+}_{\text{Profits by Strategic Bidding}} +\underbrace{\left(\Tilde{f}_i(\Tilde{g}_i^*)-f_i(\Tilde{g}^*_i)\right)}_{\text{Profits in Generation Cost}}.
    \label{eq:gain decomposition}
\end{equation}

Since for every generator $i$, $\Tilde{g}_i^*\in\{0,x,G\}$, we can define its maximal profits via strategic bidding as follows:
\begin{equation}
    \sup\{\Tilde{P}_i\}=\sup_{u\in\{0,x,G\}}\{\Tilde{P}_i|\Tilde{g}_i^*=u\}.
\end{equation}
This allows us to quantify each generation's maximal \emph{additional} profits through strategic manipulation:
\begin{equation}
    M(i)\coloneqq\ \sup\{\Tilde{P}_i\}-P_i.
    \label{def:market power}
\end{equation}
Surprisingly, $M(i)$ has a uniform and rather neat expression.

\vspace{0.1cm}
\noindent \textbf{Theorem 1}: $M(i)=c^{\{i\}}(y)-c(y)-P_i \ $.
\vspace{0.1cm}

Note that $c^{\{i\}}(y)$ denotes the minimal system cost without generator $i$'s participation. Generally, for each set $A$, we can define $c^A(y)$ as follows:
\begin{equation}
\begin{aligned}
c^{A}(y)&\coloneqq\min_{\textbf{g}} f
_k(g_k)\\
s.t. \quad &\sum\nolimits_{k=1}^n g_k=y, \\
& g_k=0, \ \ k\in A,\\
& 0\le g_k\le G, \ \ k\not\in A. 
\end{aligned}
\end{equation}
This highlights $c^{\{i\}}(y)-c(y)$ in Theorem 1 just seems like the information rent in the VCG mechanism \cite{bushnell1999international}!

\vspace{0.1cm}
\noindent \textbf{Remark:} Throughout the analysis, we assume that only one generator adopts the strategic behavior, and it knows all the information in the system. This is the standard assumption for worst case analysis and it is best suitable to understand the \emph{potential} of strategic behavior. 

\begin{figure}
    \centering
    \includegraphics[scale=0.23]{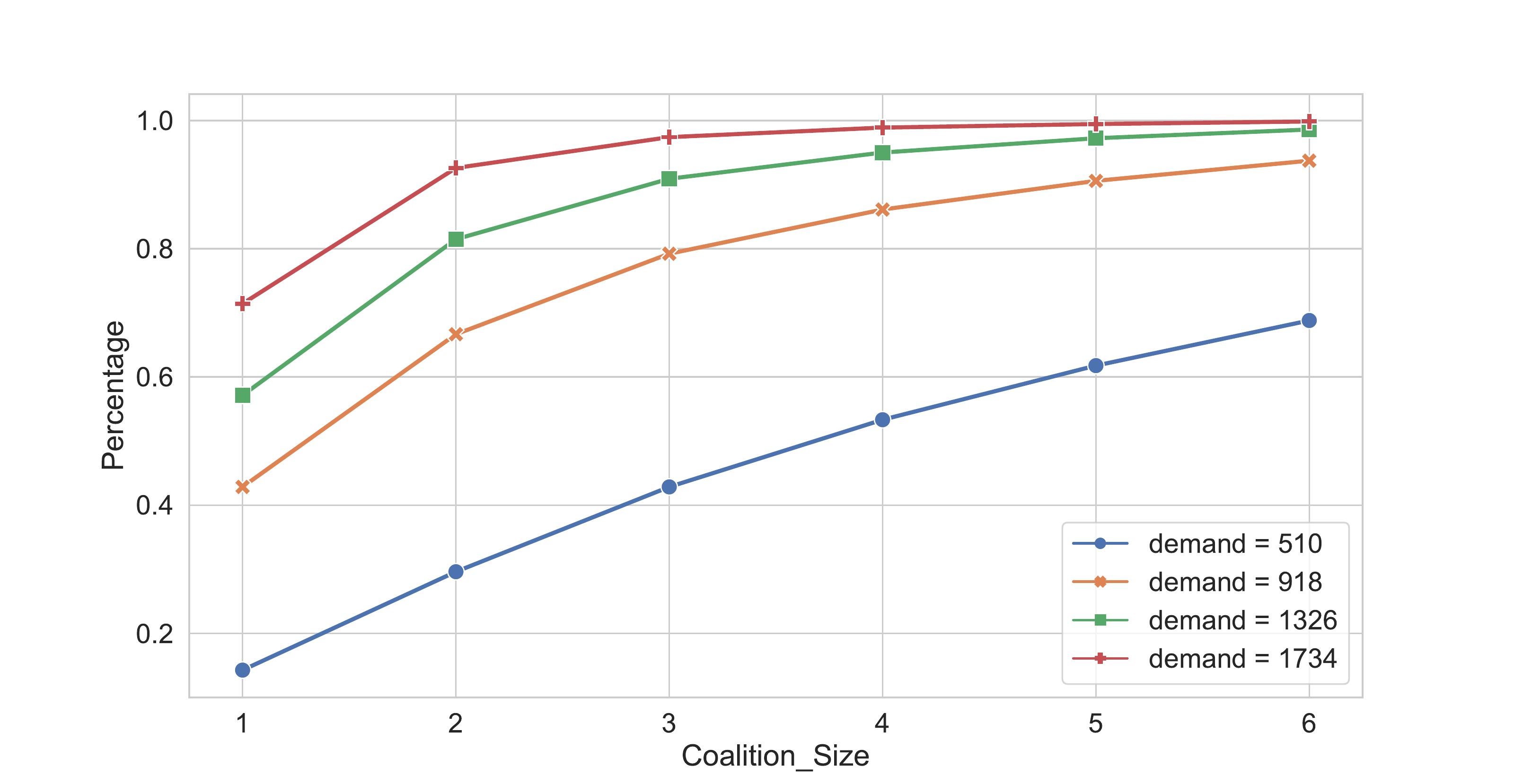}\vspace{-0.2cm}
    \caption{Percentage of coalitions with market power.\vspace{-0.5cm}}
    \label{fig:percentage}
\end{figure}

\begin{figure}
    \centering
        \includegraphics[scale=0.23]{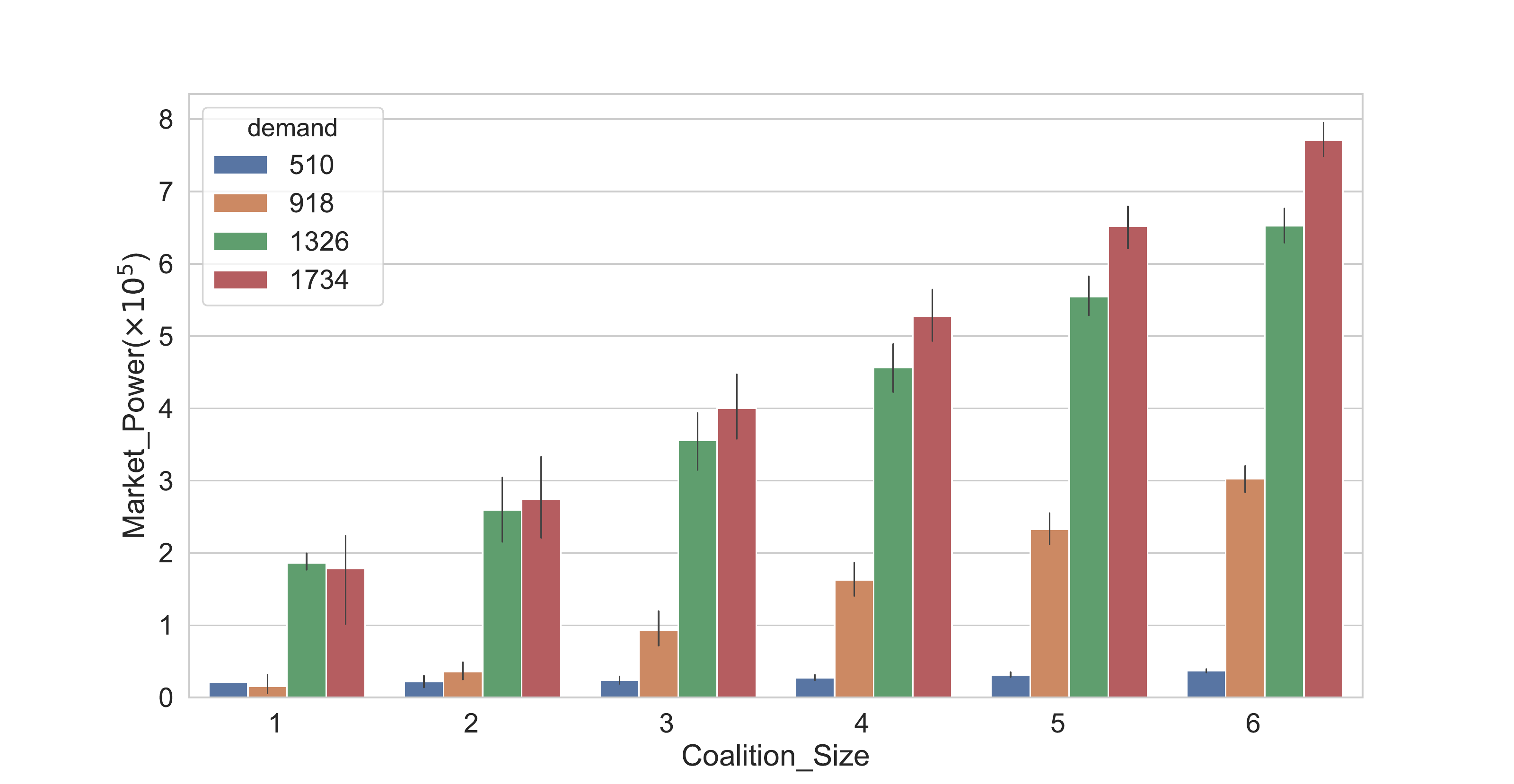}\vspace{-0.2cm}
    \caption{Market power quantification versus coalition size. \vspace{-0.5cm}}
    \label{fig:powers}
\end{figure}

\section{Market Power Index}
\label{sec:4}
The additional profits brought by strategic bidding can actually measure market power. More precisely, suppose two generators ($i$ and $j$) collude together, and all other generators are truthful, we can similarly characterize $M(i\cup j)$ as follows:
\begin{equation}
\label{M_cup_def}
   M(i\cup j)\coloneqq \sup \{\tilde{P}_{i\cup j}\}\!-\!P_{i\cup j}\!=\!c^{\{i,j\}}(y)\!-\!c(y)\!-\!P_i\!-\!P_j.
\end{equation}

Due to space limit, we omit the details in deriving Eq. (\ref{M_cup_def}). Instead, we prove the following theorem to highlight that index $M(\cdot)$ is suitable to characterize market power.

\vspace{0.1cm}
\noindent \textbf{Theorem 2}: Index $M$ is supermodular, i.e.,
\begin{equation}
    M(i\cup j)\ge M(i)+M(j), \forall i\neq j.
\end{equation}

We further examine the existence of market power in the IEEE RTS 96 node system \cite{TP-toolbox-web}. This model consists of 24 generators (which can be categorized into 7 types). We modify the system parameters such that all the generators share the same capacity of 100MW. We quantify market power with total load ranging from 510MW to 1734MW, and we consider at most 6 generators can collude together. This is because without any 6 generators, the rest can still meet the total demand. From Figs. \ref{fig:percentage} and \ref{fig:powers}, it is evident that when more generators collude together, they are more likely to have market power (Fig. \ref{fig:percentage}) and surprisingly, their average ability in manipulating the market increases almost linearly in the coalition size (Fig. \ref{fig:powers}).

\section{Concluding Remarks}
\label{sec:5}
In this letter, we identify the existence of strategic behaviors in a restricted setting for CHP scheme. This highlights the existence of market power in the general CHP schemes. We intend to examine the impacts of heterogeneous generation capacities as well as the network constraints on market power quantification in our future work.

\begin{small}
\bibliographystyle{ieeetr}
\bibliography{reference}
\end{small}

\section*{Appendix}
\label{sec:appendix}
\subsection{Proof of Lemma 1}
Given the cost structure of $f_i(g_i)$ in Eq. (\ref{eq:cost structure}), we know
\begin{equation}
    g_i^d(p)=\left\{
    \begin{aligned}
    &0 & p < f_i(G)/G\\
    &G & p\ge f_i(G)/G\\
    \end{aligned}
    \right.
\end{equation}
It is suffices to identify 
\begin{equation}
    \sum\nolimits_{i=1}^n \partial\text{Uplift}_i(p,y)/\partial p=\sum\nolimits_{i=1}^n g_i^d(p)-y.
\end{equation}
First order optimality condition yields our conclusion.

\subsection{Benchmark Profits Analysis}
\label{proof:truthful profit}
Suppose the generators are sorted according to the average cost. Then, Lemma 1 dictates 
\begin{equation}
    p^* = f_m(G)/G.
\end{equation}
Hence, for generator $i$, $i\le m$, its total profits is $f_m(G)-f_i(G)$. On the other hand, for generator $i$, $i>m$, it achieves zero profit. In summary, 
\begin{equation}
    \pi_i(p^*)=\{f_m(G)-f_i(G)\}^+.
\end{equation}

\subsection{Proof of Theorem 1}
The major difficulty lies in understanding the structure of $\sup \{\tilde{P}_i\}$. Note that, this can be analyzed case by case: $\sup \{\tilde{P}_i|\tilde{g}^*_i=0\}$, $\sup \{\tilde{P}_i|\tilde{g}^*_i=x\}$, and $\sup \{\tilde{P}_i|\tilde{g}^*_i=G\}$. Due to space limit, we only provide the full analysis for the last case. The other two cases can be analyzed in the same routine.

We can first identify that only those generators, whose rank $i$ is no greater than $m$, may have incentive to strategically bid, and the strategic bidding has to become the new determinant for the optimal price $\tilde{p}^*$. That is, after strategic bidding, generator $i$'s cost function becomes $\tilde{f}^{(m)}(G)$. This requires 
\begin{equation}
    \tilde{f}^{(m)}(G)+c^{\{i\}}(y-G)\le c^{\{i\}}(y).
\end{equation}
Hence, $\sup\{\tilde{P}_i|\tilde{g}^*_i=G\}=c^{\{i\}}(y)-c^{\{i\}}(y-G)-f_i(G)$.

Combining all the possible choices, we can show that 
\begin{equation}
    M(i)=c^{\{i\}}(y)-c(y)-P_i.
\end{equation}

\subsection{Proof of Theorem 2}

This theorem is an immediate result by combining the following two Lemmas.


\vspace{0.1cm}
\noindent\textbf{Lemma 2}: $f^{(m-1)}(G)\le c(y)-c(y-G)\le f^{(m)}(G)$.
\vspace{0.1cm}

\noindent\textbf{Proof:} We want to emphasize that this conclusion does not require truthful bidding. The key observation is to identify that given a demand of $y$, there must be one generator, whose rank is at most $m$, is dispatched to generate $G$ units. This implies $c(y)-f_m(G)\ge c(y-G)$, which constructs the second inequality.



The first inequality relies on similar observation. Given a demand of $y-G$, there must be one generator, whose rank is at least $m-1$, is dispatched to generate $0$ unit. However, when the demand is $y$, this generator can be dispatched to generate $G$. Hence, $c(y)\le f_{m-1}(G)+c(y-G)$.


\vspace{0.1cm}
\noindent \textbf{Lemma 3}: For given demand of $y$, generators $i$ and $j$, it holds
\begin{equation}
    c^{\{i,j\}}(y)-c^{\{j\}}(y)\ge c^{\{i\}}(y)-c(y).
\end{equation}

\noindent \textbf{Proof}: This lemma automatically holds when either of $g_i^*$ and $g_j^*$ is zero. Hence, without loss of generality, we assume that $g_i^*=G$, and $g_j^*=u\in\{x,G\}$. We can show that
\begin{equation}
\begin{aligned}
   \!\!\! &c^{\{i,j\}}(y)\!-\!c^{\{j\}}(y)+c(y)\!\ge\! f_{m-1}^{\{i,j\}}(G)\!-\!f_i(G)+c(y)\\
    \ge & \ c^{\{i,j\}}(y-G-u)\!+f_{m-1}^{\{i,j\}}(G)+f_j(u) \ge c^{\{i\}}(y),
\end{aligned}
\end{equation}
where $f_{m-1}^{\{i,j\}}$ denotes the $m-1^{th}$ generation cost when $i$ and $j$ do not participate. The first inequity holds due to Lemma 2, and the second inequity holds due to the following fact: 
\begin{equation}
    c(y)=c^{\{i,j\}}(y-u-G)+f_i(G)+f_j(u).
\end{equation}
The last inequity is because $c^{\{i\}}(y)$ is the optimal objective value to meet demand of $y$ without the help of generator $i$.

\end{document}